\documentclass[12pt]{amsart}
\usepackage{amscd,amssymb,epsf}
\usepackage[mathscr]{eucal}
\newtheorem{thm}{Theorem}
\newtheorem*{thm*}{Theorem}
\newtheorem{lemma}[thm]{Lemma}
\newtheorem{prop}[thm]{Proposition}
\begin{document}
\newcommand{\R}{{\mathbb R}}
\newcommand{\C}{{\mathbb C}}
\newcommand{\Id}{{\mathbb I}}
\newcommand{\Ker}{\mathsf{Ker}}
\newcommand{\Aut}{\mathsf{Aut}}
\newcommand{\Ad}{\mathsf{Ad}}
\newcommand{\Inn}{\mathsf{Inn}}
\newcommand{\tr}{\mathsf{tr}}
\renewcommand{\det}{\mathsf{det}}
\newcommand{\SL}[1]{{\mathsf{SL}}({#1})}
\newcommand{\SLt}{{\SL{2}}}
\newcommand{\GL}[1]{{\mathsf{GL}}({#1})}
\newcommand{\slt}{{\SL{2,\C}}}
\newcommand{\glt}{{\GL{2,\C}}}
\newcommand{\Mat}{{\mathsf{M}_2(\C)}}
\newcommand{\bt}{{\mathsf{t}}}
\newcommand{\hmg}{\mathsf{Hom}(\pi ,G)}
\newcommand{\W}{{\mathcal{W}}}

\title
{An exposition of results of Fricke and Vogt}
\author{William M.~Goldman}
\address{ Mathematics Department,
University of Maryland, College Park, MD  20742 USA  }
\email{ wmg@math.umd.edu }
\date{\today}
\maketitle
\tableofcontents
The purpose of this paper is to give an elementary and relatively
self-contained proof of the following basic theorem, apparently due to
Fricke~\cite{Fricke,FrickeKlein} and, earlier by Vogt~\cite{Vogt} 
Let $G=\slt$ and let $G$ act on $H=\slt\times\slt$ by
\begin{equation*}
g\cdot (\xi,\eta) \longmapsto  (g\xi g^{-1},g\eta g^{-1}).
\end{equation*}
Let $\kappa(x,y,z) = x^2 + y^2 + z^2 - xyz -2$. 

\begin{thm*}[Vogt~\cite{Vogt},Fricke~\cite{Fricke}]
Let $f:H\longrightarrow\C$ be a regular function which is invariant under
the above $G$-action. Then there exists a polynomial function 
$F(x,y,z)\in\C[x,y,z]$ such that 
\begin{equation*}
f(\xi,\eta) = F(\tr(\xi),\tr(\eta),\tr(\xi\eta)). 
\end{equation*}
For every $(x,y,z)\in\C^3$, there exists $(\xi,\eta)\in H$ such that 
\begin{equation*}
\bmatrix x \\ y \\z \endbmatrix = 
\bmatrix \tr(\xi) \\ \tr(\eta) \\\tr(\xi\eta) \endbmatrix.  
\end{equation*}
Conversely, if $\kappa(x,y,z)\neq 2$ and 
$(\xi,\eta),(\xi',\eta')\in H$ satisfy
\begin{equation}\label{eq:equalchars}
\bmatrix \tr(\xi) \\ \tr(\eta) \\ \tr(\xi\eta) \endbmatrix = 
\bmatrix \tr(\xi') \\ \tr(\eta') \\ \tr(\xi'\eta') \endbmatrix = 
\bmatrix x \\ y \\ z \endbmatrix,  
\end{equation}
then $(\xi',\eta') = g.(\xi,\eta)$ for some $g\in G$.
\end{thm*}

We regard $H$ as a {\em variety of representations\/} of a free group.
Let $\pi$ be the group freely generated by elements $X,Y$,
which we call {\em letters.\/} 
We call $X,Y,X^{-1},Y^{-1}$ {\em symbols.\/}
An element of $\pi$ is a {\em reduced word\/} $w(X,Y)$, that 
is a sequence of symbols such that no symbol follows
its inverse. The length of a reduced word $w(X,Y)$ is denoted
$\ell(w)$. The empty reduced word has length $0$ and is the identity
element $\Id\in\pi$.
Reduced words are multiplied by concatenation subject
to cancellation relations of the form 
\begin{equation*}
XX^{-1} =X^{-1}X =YY^{-1} =Y^{-1}Y  = \Id.
\end{equation*}

The Cartesian product $H$ identifies with the set $\hmg$ of all
representations $\pi\longrightarrow G$ under the correspondence
\begin{align*}
\hmg & \longrightarrow H = G \times G \\
\rho & \longmapsto (\rho(X),\rho(Y)).
\end{align*}
Since $\pi$ is generated by $X,Y$ subject to no relations,
this map is an isomorphism. Then the Vogt-Fricke theorem is the statement
that the categorical quotient of $\hmg$ by $G$ is affine 3-space $\C^3$.

This result applies to non-Euclidean geometry. Let $M$ be a compact
surface homeomorphic to a three-holed sphere.  The space of
hyperbolic- or elliptic- geometric structures on $M$ is then
parametrized by the lengths of the boundary curves, which relate to
the traces of the corresponding representation of the fundamental
group. See Buser~\cite{Buser}, Fenchel~\cite{Fenchel},
Keen~\cite{K1,K2,K3}, and Harvey~\cite{Harvey}.  In \cite{K3}, Keen
uses Fricke moduli to parametrize the Teichm\"uller space of a
one-holed torus and other surfaces.

\section{Traces in $\slt$}
\subsection{Cyclic groups}
The analogous result for one generator is the following.
\begin{thm}
Let $f:G\longrightarrow\C$ be a polynomial function invariant under
inner automorphisms of $G$. Then there exists a polynomial $F(t)\in\C[t]$
such that $f(g) = F(\tr(g))$. Conversely, if $g,g'\in G$ satisfy
$\tr(g) = \tr(g') \neq \pm 2$, then there exists $h\in G$ such that
$g' = h g h^{-1}$.
\end{thm}
\begin{proof}
Suppose $f$ is an invariant function. For $t\in\C$, define
\begin{equation*}
\xi_t := \bmatrix t & -1 \\ 1 & 0 \endbmatrix 
\end{equation*}
and define $F(t)$ by
\begin{equation*}
F(t) = f(\xi_t).
\end{equation*}
Suppose that $t\neq\pm 2$ and $\tr(g) = t$. 
Then both $g$ and $\xi_t$ have distinct eigenvalues
\begin{equation*}
\lambda_{\pm} = \frac12 \left( t \pm ( t^2 - 4)^{1/2} \right)
\end{equation*}
and $hgh^{-1} = \xi_t$ for some $h\in G$.
Thus 
\begin{equation*}
f(g) = f(h^{-1}\xi_th) = f(\xi_t) = F(t) 
\end{equation*}
as desired. If $t=\pm 2$, then by taking Jordan normal
form, either $g=\pm\Id$ or $g$ is conjugate to $\xi_t$. In the latter case,
$f(g) = F(t)$ follows from invariance. Otherwise $g$ lies in the closure
of the $G$-orbit of $\xi_t$ and $f(g) = f(\xi_t) = F(t)$ follows by
continuity of $f$.

The converse direction follows from Jordan normal form as already used above.
\end{proof}

\subsection{The Cayley-Hamilton theorem}

A principal tool in this investigation is the {\em Cayley-Hamilton theorem.\/}
If $\xi$ is a $2\times 2$-matrix,
\begin{equation}\label{eq:CH}
\xi^2 - \tr(\xi) \xi + \det(\xi)\Id = 0. 
\end{equation}

Suppose $\xi,\eta\in\slt$.
Multiplying \eqref{eq:CH} by $\xi^{-1}$ and rearranging,
\begin{equation}\label{eq:SumInvo}
\xi + \xi^{-1} = \tr(\xi) \Id 
\end{equation}
from which follows
\begin{equation}\label{eq:inverse}
\tr(\xi) =  \tr(\xi^{-1})
\end{equation}
Multiplying \eqref{eq:SumInvo} by $\eta$ and taking traces, we obtain:
\begin{thm}[The Basic Identity]
Let $\xi,\eta\in\SLt$. Then
\begin{equation}\label{eq:basic}
\tr(\xi\eta) + \tr(\xi\eta^{-1})  
= \tr(\xi)\tr(\eta).
\end{equation}
\end{thm}
\subsection{Traces of reduced words}\label{sec:wordtraces}
Here is an important special case of the Vogt-Fricke theorem.
Namely, let $w(X,Y)\in\pi$
be a reduced word. Then
\begin{align*}
H & \longrightarrow \C \\
(\xi,\eta) & \longmapsto \tr w(\xi,\eta)
\end{align*}
is a $G$-invariant function on $H$ and by the Vogt-Fricke theorem 
there exists a polynomial
$f_w(x,y,z)\in\C[x,y,z]$ such that 
\begin{equation}\label{eq:fw}
\tr w(\xi,\eta)  = f_w(\tr\xi,\tr\eta,\tr(\xi\eta))
\end{equation}
for all $\xi,\eta\in G$. We describe an algorithm for computing
$f_w(x,y,z)$. For notational convenience we shall write $t(w(X,Y))$ for
\begin{equation*}
f_{w(X,Y)}(x,y,z). 
\end{equation*}

For example, $t(\Id) = 2$ and 
\begin{align*}
t(X^{-1}) & = t(X) = x, \\
t(Y^{-1}) & = t(Y) = y, \\
t(Z^{-1}) & = t(Z) = z, \\
\end{align*}
verifying assertion \eqref{eq:fw} for words $w$ of length $\ell(w)\le 1$.
The reduced words of length two are 
\begin{equation*}
X^2,Y^2,XY,XY^{-1},YX,YX^{-1} 
\end{equation*}
 and
applications of the trace identities imply:
\begin{align*}
t(X^2) & =  x^2 - 2 \\
t(Y^2) & =  y^2 - 2 \\
t(Z^2) & =  z^2 - 2 \\
t(XY^{-1}) & =  xy - z \\
t(YZ^{-1}) & =  yz - x \\
t(ZX^{-1}) & =  zx - y 
\end{align*}
For example, for $X^2$,
\begin{equation*}
\tr(\xi^2) = \tr(\xi\xi) = \tr(\xi)\tr(\xi) - \tr(\xi\xi^{-1}) = 
\tr(\xi)^2 - 2.
\end{equation*}
and, for $XY^{-1}$,
\begin{equation*}
\tr(\xi\eta^{-1}) = \tr(\xi)\tr(\eta) - \tr(\xi\eta),
\end{equation*}
etc. For example
\begin{align}\label{eq:xyxiy}
t(XYX^{-1}Y) & = t(XY) t(X^{-1}Y) - t(X^2)  \\ 
& = z(xy - z) - (x^2 -2) \notag \\ 
&  = 2 - x^2  - z^2  + xyz. \notag
\end{align}

Assume inductively that for all reduced words $w(X,Y)\in\pi$ with
$\ell(w)<l$, 
there exists a polynomial $f_w(x,y,z) = t(w(X,Y))$ satisfying \eqref{eq:fw}.
Suppose that $u(X,Y)\in\pi$ is a reduced word of length $\ell(u) = l$.

The explicit calculations above begin the induction for $l\le 2$. 
Thus we assume $l>2$.

Furthermore, we can assume that $u$ is {\em cyclically reduced,\/} that is
the initial symbol of $u$ is not inverse to the terminal symbol of $u$.
For otherwise 
\begin{equation*}
u(X,Y) = S u'(X,Y) S^{-1}, 
\end{equation*}
where $S$ is one of the four symbols 
\begin{equation*}
X, Y,X^{-1},Y^{-1} 
\end{equation*}
and $\ell(u') = l-2$. Then $u(X,Y)$ and $u'(X,Y)$ are conjugate
and 
\begin{equation*}
t(u(X,Y)) = t(u'(X,Y)). 
\end{equation*}

If $l>2$, and $u$ is cyclically reduced,
then $u(X,Y)$ has a repeated letter, which we may assume to
equal $X$. That is, we may write 
\begin{equation*}
u(X,Y) = u_1(X,Y) u_2(X,Y) 
\end{equation*}
where $u_1$ and $u_2$ are reduced words each ending in $X^{\pm 1}$. 
Furthermore we may assume that 
\begin{equation*}
\ell(u_1) + \ell(u_2) = \ell(u) = l, 
\end{equation*}
so that $\ell(u_1) < l$ and $\ell(u_l) < l$.
Suppose first that $u_1$ and $u_2$ both end in $X$. Then 
\begin{equation*}
u(X,Y) = \left( u_1(X,Y)X^{-1}\right)X\ 
\left( u_2(X,Y)X^{-1}\right)X 
\end{equation*}
and each 
\begin{equation*}
u_1(X,Y)X^{-1}, u_2(X,Y)X^{-1} 
\end{equation*}
is reduced. Then 
\begin{equation*}
(u_1(X,Y)X^{-1})(u_2(X,Y)X^{-1})^{-1} = 
u_1(X,Y)u_2(X,Y)^{-1}
\end{equation*}
is represented by a reduced word $u_3(X,Y)$ 
satisfying $\ell(u_3) < l$.
By the induction hypothesis, there exist polynomials
\begin{equation*}
t(u_1(X,Y)), t(u_2(X,Y)), t(u_3(X,Y)) \in \C[x,y,z]
\end{equation*}
such that, for all $\xi,\eta\in G$, $i=1,2,3$,
\begin{equation*}
\tr\big( u_i(\xi,\eta) \big) = 
t\big( u_i(X,Y))(\tr(\xi),\tr(\eta),\tr(\eta) \big).
\end{equation*}
By \eqref{eq:basic}, 
\begin{equation*}
t\big(u(X,Y)\big) = t(u_1(X,Y)) t(u_2(X,Y)) - t(u_3(X,Y))
\end{equation*}
is a polynomial in $\C[x,y,z]$.
The cases when $u_1$ and $u_2$ both end in the symbols $X^{-1}, Y, Y^{-1}$ 
are completely analogous.
Since there are only four symbols, the only cyclically reduced words without
repeated symbols are commutators of the symbols, for example 
$XYX^{-1}Y^{-1}$. Repeated applications of the trace identities give:
\begin{align}\label{eq:commutator}
t(XYX^{-1}Y^{-1}) & = t(XYX^{-1}) t(Y^{-1}) - t(XYX^{-1}Y) \notag\\
& = y^2 - (2 - x^2  - z^2  + xyz) = \kappa(x,y,z). 
\end{align}
The other commutators of distinct symbols also have trace $\kappa(x,y,z)$
by identical arguments.

Thus every $w(X,Y)$ determines a polynomial $f_w(x,y,z)$ such that
\begin{equation*}
\tr w(\xi,\eta) = f_w(\tr(\xi),\tr (\eta),\tr (\xi\eta)) 
\end{equation*}
for $\xi,\eta\in\slt$.

\section{Surjectivity of characters of pairs}
We first show that
\begin{align*}
\tau: H & \longrightarrow \C^3 \\
(\xi,\eta) & 
\longmapsto \bmatrix \tr \xi \\  \tr \eta \\ \tr \xi\eta\endbmatrix. 
\end{align*}
is surjective.
Choose $\zeta\in\C$ so that
\begin{equation*}
\zeta + \zeta^{-1} = z, 
\end{equation*}
that is, $\zeta = \frac12 (z \pm \sqrt{z^2 - 4})$.
Let
\begin{equation}\label{eq:explicitrep}
\xi_x = \bmatrix x &  -1  \\ 1 & 0 \endbmatrix,\;
\eta_{(y,\zeta)} = \bmatrix 0 &  \zeta^{-1} \\ \zeta & y \endbmatrix.
\end{equation}
Then $\tau(\xi_x,\eta_{(y,\zeta)}) = (x,y,z)$.

Next we show that every $G$-invariant regular function
$f:H\longrightarrow\C$ factors through $\tau$.
To this end we need the following elementary lemma on symmetric functions:

\begin{lemma}\label{lem:even}
Let $R$ be an integral domain where $2$ is invertible,
and let $R' = R[\zeta,\zeta^{-1}]$ be the ring of Laurent
polynomials over $R$. Let $
R'\xrightarrow{\sigma} R'$ be the involution
which fixes $R'$ and interchanges $\zeta$ and $\zeta^{-1}$. Then the subring
of $\sigma$-invariants is the polynomial ring $R[\zeta+\zeta^{-1}]$.
\end{lemma}
\begin{proof}
Let $F(\zeta,\zeta^{-1})\in R[\zeta,\zeta^{-1}]$ be a $\sigma$-invariant
Laurent polynomial. 
Begin by rewriting $R'$ as the quotient of the polynomial ring
$R[x,y]$ by the ideal generated by $xy-1$. Then $\sigma$ is induced by
the involution $\tilde\sigma$ of $R[x,y]$ interchanging $x$ and $y$.
Let $f(x,y)\in R[x,y]$ be a polynomial whose image in $R'$ is $F$.
Then there exists a polynomial $g(x,y)$ such that 
\begin{equation*}
f(x,y) - f(y,x) = g(x,y) (xy - 1). 
\end{equation*}
Clearly $g(x,y) = - g(y,x)$. 
Let 
\begin{equation*}
\tilde f(x,y) = f(x,y) - \frac12\; g(x,y) (x y - 1)
\end{equation*}
so that $\tilde f(x,y)  = \tilde f(y,x)$. By the theorem on elementary
symmetric functions,
\begin{equation*}
\tilde f(x,y) = h(x + y, xy) 
\end{equation*}
for some polynomial $h(u,v)$.
Therefore $F(\zeta,\zeta^{-1}) = h(\zeta + \zeta^{-1},1)$ as desired.
\end{proof}

By definition $f(\xi,\eta)$ is a polynomial in the matrix entries of
$\xi$ and $\eta$; two polynomials which differ by elements in the
ideal generated by $\det(\xi)-1$ and $\det(\eta)-1$ are regarded as
equal. Thus $f(\xi_x,\eta_{(y,\zeta)})$ equals a function $g(x,y,\zeta)$
which is a polynomial in $x,y\in\C$ and a Laurent polynomial in
$\zeta\in\C^*$.

\begin{lemma}\label{lem:commonperp}
Let $\xi,\eta\in G$ such that $\kappa(\tau(\xi,\eta))\neq 2$. 
Then there exists
$g\in G$ such that 
\begin{equation*}
g\cdot (\xi,\eta) =  (\xi^{-1},\eta^{-1}).
\end{equation*}
\end{lemma}
\begin{proof}
Let $(x,y,z) = \tau(\xi,\eta)$. By a simple calculation
\begin{equation*}
\tr [\xi,\eta] = \kappa(x,y,z) 
\end{equation*}
where $[\xi,\eta] = \xi\eta\xi^{-1}\eta^{-1}$.

Let $L = \xi\eta-\eta\xi$. 
(Compare \S 4 of J\o rgensen~\cite{J} or Fenchel~\cite{Fenchel}.)
Then $\tr(L) = \tr(\xi\eta) -\tr(\eta\xi) = 0$. Furthermore
for any $2\times 2$ matrix $M$, the characteristic polynomial
\begin{equation*}
\lambda_M(t) := \det(t\Id - M) = t^2 - \tr(M) t + \det(M).
\end{equation*}
Thus
\begin{align*}
\det(L) & = \det( [\xi,\eta] - I)\det(\eta\xi) \\ & = \det([\xi,\eta]-I) \\ 
& = - \lambda_{[\xi,\eta]}(1) \\ 
& = - 2 + \tr [\xi,\eta] \\ 
& = - 2 + \kappa(x,y,z) \neq 0.
\end{align*}
Choose $\mu\in\C^*$ such that $\mu^2 \det(L) = 1$ and let $g = \mu L\in G$.

Since $\tr(g) = 0$ and $\det(g) =1$, the Cayley-Hamilton Theorem 
$\lambda_M(M) = 0$ implies that $g^2 = -\Id$. Similarly
\begin{align*}
\det(g\xi) &= \det(g) = 1 \\ 
\det(g\eta) &= \det(g) = 1, 
\end{align*}
and 
\begin{align*}
\tr(g\xi) & = \mu (\tr((\xi\eta)\xi) - \tr((\eta\xi)\xi)) \\
& = \mu (\tr(\xi(\eta\xi)) - \tr((\eta\xi)\xi)) = 0 
\end{align*}
and 
\begin{align*}
\tr(g\eta) & = \mu (\tr((\xi\eta)\eta) - \tr((\eta\xi)\eta)) \\ 
& = \mu (\tr((\xi\eta)\eta) - \tr(\eta(\xi\eta)) = 0
\end{align*}
so $(g\xi)^2 = (g\eta)^2 = -\Id$. Thus
$ g \xi g^{-1} \xi = 
- g\xi g\xi = \Id$
whence $g \xi g^{-1} = \xi^{-1}.$ 
Similarly $g \eta g^{-1} = \eta^{-1}$, 
concluding the proof of the lemma. \end{proof}

Apply Lemma~\ref{lem:commonperp} 
to $\xi = \xi_x$ and $\eta = \eta_{(y,\zeta)}$ as above to obtain
$g$ such that conjugation by $g$ maps
\begin{equation*}
\xi  \longmapsto
\xi^{-1} =  \bmatrix 0 &  1  \\ -1 & x \endbmatrix
\end{equation*}
and 
\begin{equation*}
\eta  \longmapsto
\eta^{-1} =  \bmatrix y &  - 1/\zeta  \\ \zeta & 0 \endbmatrix.
\end{equation*}
Let 
\begin{equation*}
h = \bmatrix 0 & 1 \\ 1 & 0 \endbmatrix 
\end{equation*}
Then 
\begin{equation*}
hg \xi (hg)^{-1} = h \xi^{-1} h^{-1} =  
\bmatrix x &  -1  \\ 1 & 0 \endbmatrix = \xi
\end{equation*}
and 
\begin{equation*}
hg \eta (hg)^{-1} = 
h \eta^{-1} h^{-1} =  \bmatrix 0 &  \zeta  \\ -1/\zeta & y \endbmatrix.
\end{equation*}
Thus
\begin{align*}
g(x,y,\zeta) & =  f(\xi,\eta) \\ & 
= f(hg \xi (hg)^{-1},hg \eta (hg)^{-1}) \\ 
& =  g(x,y,\zeta^{-1}).
\end{align*}
Lemma~\ref{lem:even} 
implies that there exists a polynomial $F\in\C[x,y,z]$ such that
\begin{equation}\label{eq:even}
g(x,y,\zeta) = F(x,y,\zeta + 1/\zeta)
\end{equation}
whenever $\kappa(x,y,\zeta+1/\zeta)\neq 2$.
Since this condition defines a nonempty Zariski-dense open set,
\eqref{eq:even} holds on all of $\C^2\times\C^*$ and
$f(\xi,\eta) = F(\tr(\xi),\tr(\eta), \tr(\xi\eta))$ as claimed.

\section{Injectivity of characters of pairs}

Finally we show that if $(\xi,\eta),(\xi',\eta')\in H$ 
satisfy \eqref{eq:equalchars}
and $\kappa(x,y,z)\neq 2$, then $(\xi,\eta)$ and $(\xi',\eta')$ 
are $G$-equivalent. 
By \S\ref{sec:wordtraces}, the triple 
\begin{equation*}
\bmatrix x \\ y \\ z\endbmatrix  = 
\bmatrix \tr \xi \\ \tr \eta \\ \tr (\xi\eta)\endbmatrix 
\end{equation*}
determines the character function
\begin{align*}
\chi:\pi &\longrightarrow \C  \\
w(X,Y) & \longmapsto \tr w(\xi,\eta) = f_w(x,y,z).
\end{align*}
Let $\rho$ and $\rho'$ denote the representations $\pi\longrightarrow G$
taking $X,Y$  to $\xi,\eta$ and $\xi',\eta'$ respectively
and let $\chi,\chi'$ denote their respective characters. 
Then our hypothesis 
\eqref{eq:equalchars} implies that $\chi=\chi'$.

\begin{lemma}\label{lem:kappaIrr}
Let $\rho:\pi\longrightarrow\slt$ be a representation and
$(x,y,z)$ be as above. Then $\rho$ is irreducible
if and only if $\kappa(x,y,z)\neq 2$.
\end{lemma}
\begin{proof}
Suppose first that $\rho$ is reducible.
If $X,Y$ generate a representation
with an invariant subspace of $\C^2$ of dimension one, this representation
is conjugate to one in which $\rho(X)$ and $\rho(Y)$ are 
upper-triangular. 
Denoting their diagonal entries by $a,a^{-1}$ and $b,b^{-1}$ respectively,
the diagonal entries of $\rho(XY)$ are 
$ab^{\pm 1}, a^{-1}b^{\mp 1}$. 
Thus
\begin{align*}
x & = a + a^{-1}, \\ y & = b  + b^{-1}, \\
z & = a b^{\pm 1} + a^{-1} b^{\mp 1}.
\end{align*}
By direct computation, $\kappa(x,y,z) = 2$. 

Conversely, suppose that $\kappa(x,y,z) = 2$. 
Let $\mathfrak{A}\subset M_2(\C)$ denote the linear span of
$\Id,\xi,\eta,\xi\eta$. 
Identities derived from the Cayley-Halmilton theorem
\eqref{eq:CH} such as \eqref{eq:SumInvo} imply that
$\mathfrak{A}$ is a subalgebra of $M_2(\C)$. 
For example, $\xi^2$ is the linear combination $ -1 + x \xi$ and
\begin{equation*}
\eta\xi = (z - xy) 1 +  y \xi - x \eta - \xi\eta.
\end{equation*}

In the basis of $M_2(\C)$ by elementary matrices, the map
\begin{align*}
\C^4 & \longrightarrow M_2(\C) \\
\bmatrix x_1 \\ x_2 \\ x_3 \\ x_4 \endbmatrix & \longmapsto
x_1 \Id + x_2 \xi + x_3 \eta + x_4 \xi\eta
\end{align*}
has determinant $2-\kappa(x,y,z) = 0$ and is not surjective. 
Thus $\mathfrak{A}$ is a proper subalgebra of $M_2(\C)$ and 
the representation is reducible.
\end{proof}

Thus $\rho$ and $\rho'$ are irreducible representations on $\C^2$.
Burnside's Theorem (see Lang~\cite{Lang}, p.445) implies 
the corresponding representations 
(also denoted $\rho, \rho'$ respectively)
of the group algebra $\C\pi$ into $M_2(\C)$ are surjective. 
Since the trace form
\begin{align*}
M_2(\C) \times  M_2(\C) &\longrightarrow \C \\
(A,B) & \longmapsto \tr(AB)
\end{align*}
is nondegenerate, the kernel $K$ of $\rho:\C\pi\longrightarrow M_2(\C)$
consists of all 
\begin{equation*}
\sum_{\alpha\in\pi} a_\alpha \alpha \in \C\pi
\end{equation*}
such that 
\begin{align*}
0 & = \tr \Bigg( \Big(\sum_{\alpha\in\pi} a_\alpha \rho(\alpha)\Big) 
\rho(\beta) \Bigg) \\ 
& = \sum_{\alpha\in\pi} a_\alpha \tr\rho(\alpha\beta) \\
& = \sum_{\alpha\in\pi} a_\alpha \chi_{\alpha\beta}
\end{align*}
for all $\beta\in\pi$.
Thus the kernels of both representations of $\C\pi$ are equal, and
$\rho$ and $\rho'$ induce algebra isomorphisms
\begin{equation*}
\tilde\rho, \tilde\rho':\C\pi/K \longrightarrow M_2(\C)
\end{equation*}
respectively. The composition $\tilde\rho'\circ\rho^{-1}$ is an
automorphism of the algebra $M_2(\C)$, which must be induced by
conjugation by $g\in\glt$. (See, for example, Corollary 9.122, p.734
of Rotman~\cite{Rotman}.) In particular 
$\rho'(\gamma) = g \rho(\gamma)g^{-1}$
as desired.

\section{Trace relations for triples}
Finally, we consider the case of triples $(\xi,\eta,\zeta)\in G$,
or equivalently representations of the free group of rank three.
Unlike the rank two case, the quotient is no longer an affine space.
Rather, the coordinate ring of the quotient has dimension six, but generated
by the eight functions
\begin{equation*}
f_X, f_Y, f_Z, f_{XY}, f_{YZ}, f_{ZX}, f_{XYZ}, f_{XZY} 
\end{equation*}
subject to the two relations expressing the sum and product of traces
of the length 3 monomials in terms of traces of monomials of length 1
and 2:

\begin{align*}
f_{XYZ} + f_{XZY} & =  f_{XY} f_Z + f_{YZ} f_X +f_{ZX} f_Y + f_X f_Y f_Z \\
f_{XYZ} \;  f_{XZY}  & =   (f_X)^2 + (f_Y)^2 + (f_Z)^2   \\
& \quad + f_{XY}^2 + f_{YZ}^2 + f_{ZX}^2 \notag \\ 
& \quad \quad -  (
f_Xf_Yf_{XY} + f_Yf_Zf_{YZ} + f_Zf_Xf_{ZX})
\notag \\ 
&  \quad \quad \quad + f_{XY}f_{YZ}f_{ZX} - 4 \notag
\end{align*}
We call these two identities Fricke's Sum and Product Relations respectively.

\subsection*{Fricke's Sum Relation}
\begin{equation}\label{eq:sum}
f_{XYZ} + f_{XZY}  =  f_{XY} f_Z + f_{YZ} f_X +f_{ZX} f_Y + f_X f_Y f_Z
\end{equation}

To prove this formula, apply the Basic Identity three times:
\begin{align}
t(XYZ) + t(XYZ^{-1}) & = t(XY)t(Z) \label{eq:xyz}\\ 
t(Z^{-1}XY) + t(Z^{-1}XY^{-1}) & = t(Z^{-1}X)t(Y) 
				\label{eq:zxy}\\ 
& = (t(Z)t(X)-t(ZX))t(Y)  \notag \\
t(Y^{-1}Z^{-1}X) + t(Y^{-1}Z^{-1}X^{-1}) & = 
		t(Y^{-1}Z^{-1}) t(X) \label{eq:yzx}\\
& = t(YZ) t(X)\notag\end{align}
Now subtract  \eqref{eq:zxy} from the sum of \eqref{eq:xyz}
and \eqref{eq:yzx}, using the facts
\begin{align*}
t(XYZ^{-1}) & = t(Z^{-1}XY) \\
t(Z^{-1}XY^{-1}) & = t(Y^{-1}Z^{-1}X) \\
t(Z^{-1}X^{-1}Y^{-1}) & = 
t(Y^{-1}Z^{-1}X^{-1}) = 
t(XZY) 
\end{align*}
to obtain
\begin{align*}
t(XYZ) + t(XZY) & = t(XY)t(Z) \\ & \qquad - t(Z)t(X)t(Y) + t(ZX)t(Y)\\ 
& \qquad \qquad + t(YZ)t(X)
\end{align*}
from which \eqref{eq:sum} follows.


\subsection*{Fricke's Product Relation}
\begin{align}\label{eq:product}
f_{XYZ} \;  f_{XZY}  & =   (f_X)^2 + (f_Y)^2 + (f_Z)^2   \\
& \qquad + f_{XY}^2 + f_{YZ}^2 + f_{ZX}^2 \notag \\ 
& \qquad \qquad -  (
f_Xf_Yf_{XY} + f_Yf_Zf_{YZ} + f_Zf_Xf_{ZX})
\notag \\ 
&  \qquad \qquad \qquad + f_{XY}f_{YZ}f_{ZX} - 4 \notag
\end{align}

We derive this formula in several steps. A direct application of the
Basic Identity \eqref{eq:basic} yields:
\begin{align}\label{eq:zxzy}
t(ZXZY) & = t(ZX) t(ZY) - t(XY^{-1}) \\
& = t(ZX)t(ZY) - \left( t(X)t(Y) - t(XY) \right) \notag\\
& = t(ZX)t(ZY) - t(X)t(Y) + t(XY) \notag 
\end{align}
By \eqref{eq:xyxiy} applied to $X,Z^{-1}$:
\begin{equation}\label{eq:xzixizi}
t(XZ^{-1}X^{-1}Z^{-1})  =  
t(X)t(Z)t(ZX) - t(ZX)^2 - t(X)^2 + 2
\end{equation}


\begin{align}\label{eq:xyzxzy}
t(XYZXZY) & = t(XY) t(ZXZY)  \\ & \qquad - t(XZ^{-1}X^{-1}Z^{-1}) \notag\\
& = t(XY) \left( t(ZX)t(ZY) - t(X)t(Y) + t(XY) \right) \notag \\
& \qquad - \left(  t(X)t(Z)t(ZX) - t(ZX)^2 - t(X)^2 + 2 \right) \notag\\
& \qquad\qquad\qquad\text{(by \eqref{eq:zxzy} and \eqref{eq:xzixizi})}\notag\\
%
& = t(XY)t(ZX)t(YZ) \notag \\ 
& \qquad \quad - t(X)t(Y)t(XY) - t(Z)t(X)t(ZX) \notag \\ 
& \qquad \qquad + t(XY)^2 + t(X)^2 - 2  \notag
\end{align}

\noindent Finally, appplying \eqref{eq:xyzxzy} and the Commutator Identity
\eqref{eq:commutator} to $Y,Z$:
\begin{align*}
t(XYZ)t(XZY) & = t(XYZXZY) + t(YZY^{-1}Z^{-1}) \\
& = 
\Big( t(XY)t(ZX)t(YZ)  - t(X)t(Y)t(XY) \\ 
& \quad  - t(Z)t(X)t(ZX) + t(XY)^2 + t(ZX)^2 + t(X)^2 - 2 \Big) \\
& \qquad + \Big(
t(Y)^2 + t(Z)^2 + t(YZ)^2 - t(Y)t(Z)t(YZ) - 2  \Big) \end{align*}
from which \eqref{eq:product} follows.

\subsection{The coordinate ring is a quadratic extension}
If $\pi$ is freely generated by $X,Y,Z$, then
the ring of $G$-invariant polynomials on $\hmg$ is a quadratic
extension of the polynomial ring 
\begin{equation*}
\C[f_X,f_Y,f_Z,f_{XY},f_{YZ},f_{ZX}]. 
\end{equation*}
The algebraic generator $\lambda$ (which corresponds to
$f_{XYZ}$ or $f_{XZY}$) satisfies the quadratic equation
\begin{equation*}
\lambda^2 - f_\Sigma \lambda + f_\Pi = 0
\end{equation*}
where
\begin{equation*}
f_\Sigma = f_{XY} f_Z + f_{YZ} f_X +f_{ZX} f_Y - f_X f_Y f_Z 
\end{equation*}
and 
\begin{align*}
f_\Pi & =  (f_X^2 + f_Y^2 + f_Z^2) \\
& \quad + (f_{XY}^2 + f_{YZ}^2 + f_{ZX}^2) \\
& \qquad - (f_Xf_Yf_{XY} + f_Yf_Zf_{YZ} + f_Zf_Xf_{ZX}) \\
& \qquad \quad + f_{XY}f_{YZ}f_{ZX} - 4 
\end{align*}
(Compare Magnus~\cite{Magnus}.)

\section{Surjectivity in rank 3}

Let V be the $\SLt$-character variety of the (rank three) free group
$\pi=\langle A_1,A_2,A_3\rangle$. In trace coordinates
\begin{equation*}
t_I([\rho]) := \tr \big(
\rho(A_{i_1} A_{i_2} \dots A_{i_k} )\big)
\end{equation*}
where
\begin{equation*}
I = (i_1, i_2, \dots i_k ),
\end{equation*}
$V\subset\C^8$ is the codimension two subvariety 
defined by the two equations:
\begin{align}
t_{123} \,+\, t_{132} & \;=\; t_{12}t_3 + t_{13}t_2 + t_{23}t_1 + t_1t_2t_3 
\label{eq:sum1}  \\
t_{123}\,\,\,  t_{132} & \;=\;  (t_1^2 + t_2^2 +t_3^2) \,+\,
(t_{12}^2 + t_{23}^2 \,+\, t_{13}^2) \  +  \label{eq:product1}\\  
& \qquad \  ( t_1t_2t_{12} + t_2t_3t_{23} + t_3t_1t_{13} ) \,+\,
t_{12} t_{23} t_{13} \,-\, 4. \notag 
\end{align}
By eliminating $t_{132}$ in \eqref{eq:sum1} as
\begin{equation*}
 t_{132}  \;=\; t_{12}t_3 + t_{13}t_2 + t_{23}t_1 + t_1t_2t_3  - t_{123},
\end{equation*}
the variety $V$ may be expressed as the hypersurface in $\C^7$
consisting of all
\begin{equation*}
\big(t_1,t_2,t_3,t_{12},t_{23},t_{13}\big) \in \C^7
\end{equation*}
satisfying
\begin{align*}
t_{123} \;
\big( t_{12}t_3 + t_{13}t_2 + t_{23}t_1 + t_1t_2t_3  - t_{123}\big) & = \\
  (t_1^2 + t_2^2 +t_3^2) \,+\,
(t_{12}^2 + t_{23}^2 \,+\, t_{13}^2) \  +  & 
 ( t_1t_2t_{12} + t_2t_3t_{23} + t_3t_1t_{13} ) \,+\,
t_{12} t_{23} t_{13} \,-\, 4. 
\end{align*}

\begin{prop}\label{prop:onto}
The projection
\begin{align*}
V & \xrightarrow{\bt}  \C^6 \\
[\rho] &\longmapsto \bmatrix 
t_1(\rho) \\ t_2(\rho) \\ t_3(\rho) \\ 
t_{12}(\rho) \\ t_{23}(\rho) \\ t_{13}(\rho)
\endbmatrix
\end{align*}
is surjective.
\end{prop}

The polynomial ring  
\begin{equation*}
 \C[t] := \C[t_1,t_2,t_3,t_{12},t_{23},t_{13}]
\end{equation*} 
is the coordinate ring of affine space $\C^6$.
Projection $V \xrightarrow{\bt}\C^6$ 
is a 2-1 map corresponding to the quadratic extension
\begin{equation*}
\C[V]\; \cong \; \C[t][z]\, \bigg/\, \bigg(z^2 - P(t) z + Q(t) \bigg)
\end{equation*}
where $P,Q$ are the polynomials appearing above:
\begin{align*}
P(t) & = t_{12}t_3 + t_{13}t_2 + t_{23}t_1 + t_1t_2t_3  \\
Q(t) & = (t_1^2 + t_2^2 +t_3^2) +
(t_{12}^2 + t_{23}^2 +t_{13}^2) + \\
& \qquad( t_1t_2t_{12} + t_2t_3t_{23} + t_3t_1t_{13} ) +
t_{12} t_{23} t_{13} - 4.
\end{align*}
\begin{proof}
By the Vogt-Fricke theorem in rank two, 
there exist $A_1,A_2\in\SLt$ such that
\begin{align}\label{eq:initial}
\tr(A_1) & = t_1, \notag \\ 
\tr(A_2) &= t_2,   \\ 
\tr(A_1A_2) &= t_{12}\notag.
\end{align}
We seek $A_3\in\SLt$ such that
\begin{align}\label{eq:rankthreeequations}
\tr(A_3) & = t_3, \notag\\  
\tr(A_2A_3) &= t_{23},\notag \\  
\tr(A_1A_3) &= t_{13}. 
\end{align}
To this end, consider the affine subspace $\W$ of $\Mat$ consisting
of matrices $W$ satisfying
\begin{align}\label{eq:threeconditions}
\tr(W) & = t_3, \notag\\ 
\tr(A_2W) & = t_{23}, \\  
\tr(A_1W) & = t_{13}.\notag   
\end{align}
Since the bilinear pairing 
\begin{align*}
\Mat \times \Mat & \longrightarrow \C \\
(X,Y) & \longmapsto \tr(XY)
\end{align*}
is nondegenerate, each of the three equations in \eqref{eq:threeconditions}
describes an affine hyperplane in $\Mat$.
We first suppose that $\kappa(t_1,t_2,t_{12})\neq 2$,
that is, 
\begin{equation}\label{eq:irreducible}
4 - t_1^2 -t_2^2 - t_{12}^2 + t_1t_2t_{12} \neq 0. 
\end{equation}

\begin{lemma}
There exists $A_1,A_2\in\SLt$ satisfying \eqref{eq:initial} 
such that $\W$ is nonempty.
\end{lemma}
\begin{proof}
Since $\kappa(t_1,t_2,t_{12})\neq 2$. 
the pair $A_1,A_2$ generates an irreducible representation. 
In particular $\{\Id,A_1,A_2\}$ is a 
linearly independent subset of the $4$-dimensional vector space $\Mat$.
The three conditions of \eqref{eq:threeconditions} are independent,
so $\W\subset\Mat$ is an affine line.
\end{proof}
Let $W_0,W_1\in\W$ be distinct elements in this line. Then the function
function
\begin{align*}
\C &\longrightarrow \C \\
 s &\longmapsto \det \big(s W_1 + (1-s) W_0\big) 
\end{align*}
is polynomial of degree $\le 2$, and is thus onto unless it is constant. 

\begin{lemma} Let $W_0,W_1\in\Mat$. Then
\begin{align*}
\det\big(W_0 & + s (W_1-W_0)\big)  = \\ 
& \det(W_0)  \,+\,  s\; \big( \tr(W_0)\tr(W_1-W_0)-\tr\big(W_0(W_1-W_0)\big) \\
& \qquad\qquad  +\,s^2\; \det(W_1-W_0)
\end{align*}
\end{lemma}
\begin{proof}
Clearly
\begin{equation}\label{eq:trlinear}
\tr\big(W_0  \,+\; s\, (W_1-W_0)\big)  =
\tr(W_0) \,+\; s \, \tr(W_1-W_0).
\end{equation}
Now
\begin{equation}\label{eq:detmat}
\det(W) = \frac{\tr(W)^2-\tr(W^2)}2 
\end{equation}
whenever $W\in\Mat$.
Now apply \eqref{eq:detmat} to \eqref{eq:trlinear}
taking 
\begin{equation*}
W = W_0 + s (W_1-W_0) 
\end{equation*}
\end{proof} 

Thus the restriction $\det|_\W$ is constant only if $\det(W_1-W_0) = 0$.

Work in the slice
\begin{equation*}
A_1 = \bmatrix t_1 & -1 \\ 1 & 0 \endbmatrix,\;
A_2 = \bmatrix 0 & \xi \\ -\xi^{-1} & t_2 \endbmatrix,
\end{equation*}
where $\xi + \xi^{-1} = t_{12}$. 

The matrix $W_0\in\SLt$ defined by:
\begin{equation*}
W_0  = 
\bmatrix t_3 & 
\big( (t_{13}-t_1t_3)\xi + t_{23}\big)\xi/(\xi^2-1) \\
\big( (t_{13}-t_1t_3) + t_{23}\xi\big)/(\xi^2-1) & 0 
\endbmatrix, 
\end{equation*}
satisfies  \eqref{eq:threeconditions}.
Any other $W\in\W$ must satisfy
\begin{align}\label{eq:threelinearconditions}
\tr\big(W-W_0\big) & = 0, \notag\\ 
\tr\big(A_2(W-W_0)\big) & = 0, \\  
\tr\big(A_1(W-W_0)\big) & = 0.\notag   
\end{align}
\begin{lemma} 
Any solution $W-W_0$ of \eqref{eq:threelinearconditions} is a multiple
of
\begin{equation*}
[A_1,A_2] =
\bmatrix \xi^{-1} - \xi & -t_2 + t_1\xi \\
-t_2 + \xi^{-1}t_1 & 
 \xi  - \xi^{-1}  \endbmatrix.
\end{equation*}
\end{lemma} 
\begin{proof} 
The first equation in \eqref{eq:threelinearconditions} 
asserts that $W-W_0$ lies in the subspace $\slt$, upon which the
trace form is nondegenerate. The second and third equations assert
that $W-W_0$ is orthogonal to $A_1$ and $A_2$. By \eqref{eq:irreducible},
$A_1,A_2$ and $\Id$ are linearly independent in $\Mat$, so the solutions
of \eqref{eq:threelinearconditions} form a one-dimensional linear subspace.
The Lie product
\begin{equation*}
[A_1,A_2] = A_1 A_2 - A_2 A_1  
\end{equation*}
is nonzero, lies in $\slt$, and since
\begin{equation*}
\tr( A_i A_1 A_2) = \tr( A_i A_2 A_1)
\end{equation*}
for $i=1,2$, orthogonal to $A_1$ and $A_2$. The lemma follows.
\end{proof} 
Parametrize $\W$  explicitly as
\begin{equation*}
W = W_0 + s [A_1,A_2].
\end{equation*}
Since
\begin{equation*}
\det([A_1,A_2) = 4 - (t_1^2 +t_2^2 +t_{12}^2 - t_1t_2t_{12}) =
2 - \kappa(t_1,t_2,t_{12}) \neq 0, 
\end{equation*}
Thus the polynomial
\begin{equation*}
\W \xrightarrow{\det} \C 
\end{equation*}
is nonconstant, and hence onto.
Taking $W_1\in\det_\W^{-1}(1)$, the proof is complete assuming
\eqref{eq:irreducible}.


Now consider the case when $4 - t_1^2 -t_2^2 - t_{12}^2 + t_1t_2t_{12} = 0$.
In that case there exist $a_1,a_2\in\C^*$ such that
\begin{equation*}
t_i = a_i + (a_i)^{-1} 
\end{equation*}
for $i=1,2$. Then either
\begin{equation}\label{eq:bothplus}
t_{12} = a_1a_2 + (a_1a_2)^{-1} 
\end{equation}
or
\begin{equation}\label{eq:plusminus}
t_{12} = a_1(a_2)^{-1} + (a_1)^{-1}a_2. 
\end{equation}
In the first case \eqref{eq:bothplus}, set
\begin{align*}
A_1 &:=  \bmatrix a_1 & t_{13} - a_1 t_3 \\ 0 & (a_1)^{-1} \endbmatrix, \\
A_2 &:=  \bmatrix a_2 & t_{23} - a_2 t_3  \\ 
                   0  & (a_2)^{-1} \endbmatrix, \\
A_3 &:=  \bmatrix t_3 & -1 \\ 1 & 0 \endbmatrix
\end{align*}
and in the second case \eqref{eq:bothplus}, set
\begin{align*}
A_1 &:=  \bmatrix (a_1)^{-1} & t_{13} - (a_1)^{-1} t_3 
\\ 0 & a_1 \endbmatrix, \\
A_2 &:=  \bmatrix a_2 & t_{23} - a_2 t_3  \\ 
                   0  & (a_2)^{-1} \endbmatrix, \\
A_3 &:=  \bmatrix t_3 & -1 \\ 1 & 0 \endbmatrix
\end{align*}
obtaining $(A_1,A_2,A_3)\in\SLt$ 
explicitly solving \eqref{eq:initial} and \eqref{eq:rankthreeequations}.
This completes the proof of Proposition~\ref{prop:onto}..
\end{proof}

%
%



Writing
\begin{equation*}
W = a_0 I + a_1 A_1 + a_2 A_2 
\end{equation*}

\subsubsection*{Acknowlegement}
I am grateful to Hyman Bass, Richard Brown, Carlos Florentino, Linda
Keen and Sean Lawton for their interest in this manuscript, and for
providing valuable suggestions.

%
%
%
%
%
%
%
%
%
\makeatletter \renewcommand{\@biblabel}[1]{\hfill#1.}\makeatother
\renewcommand{\bysame}{\leavevmode\hbox to3em{\hrulefill}\,}

\end{document}